\documentclass[12pt,reqno,british]{amsart}

\usepackage[utf8]{inputenc}
\usepackage[T1]{fontenc}
\usepackage{babel,mathtools,amssymb,amscd,url,paralist,xspace,mathrsfs,textcmds,hyperref,microtype,stmaryrd}

\mathtoolsset{mathic=true}
\setdefaultenum{\upshape (i)}{}{}{}

\theoremstyle{plain}
\newtheorem{theorem}{Theorem}[section]
\newtheorem{proposition}[theorem]{Proposition}
\newtheorem{lemma}[theorem]{Lemma}
\newtheorem{corollary}[theorem]{Corollary}

\theoremstyle{definition}
\newtheorem{definition}[theorem]{Definition}

\theoremstyle{remark}
\newtheorem{Example}[theorem]{Example}

\newtheorem*{acknowledgements}{Acknowledgements}

\newenvironment{example}{\begin{Example}\pushQED{\qee}}{\popQED\end{Example}}
\DeclareRobustCommand{\qee}{%
  \ifmmode \mathqee
  \else
    \leavevmode\unskip\penalty9999 \hbox{}\nobreak\hfill
    \quad\hbox{\qeesymbol}%
  \fi
}
\newcommand{\mathqee}{\quad\hbox{\qeesymbol}}
\newcommand{\qeesymbol}{\ensuremath\diamondsuit}

\numberwithin{equation}{section}
\numberwithin{figure}{section}

\newcommand{\itref}[1]{\eqref{#1}}

\newcommand{\Lie}[1]{\operatorname{\textsl{#1}}}
\newcommand{\lie}[1]{\operatorname{\mathfrak{#1}}}
\newcommand{\GL}{\Lie{GL}}
\newcommand{\SL}{\Lie{SL}}

\newcommand{\SU}{\Lie{SU}}

\newcommand{\An}{A^n}
\newcommand{\an}{\lie a_n}
\newcommand{\AM}{A_M}
\newcommand{\aM}{\lie a_M}
\newcommand{\AP}{A_P}
\newcommand{\aP}{\lie a_P}
\newcommand{\AW}{A_W}
\newcommand{\aW}{\lie a_W}
\newcommand{\Tn}[1][n]{T^{#1}}
\newcommand{\tn}{\lie t_n}
\newcommand{\TM}{T_M}
\newcommand{\tM}{\lie t_M}

\newcommand{\tP}{\lie t_P}
\newcommand{\XM}{\lie X}

\newcommand{\Ln}{\Lambda_n}
\newcommand{\LM}{\Lambda_M}
\newcommand{\LP}{\Lambda_P}

\newcommand{\bC}{\mathbb C}
\newcommand{\bH}{\mathbb H}
\newcommand{\bQ}{\mathbb Q}
\newcommand{\bR}{\mathbb R}
\newcommand{\bZ}{\mathbb Z}

\newcommand{\sA}{\mathscr A}
\newcommand{\cF}{\mathcal F}
\newcommand{\cL}{\mathcal L}

\DeclareMathOperator{\colim}{colim}
\DeclareMathOperator{\CP}{\mathbb CP}
\DeclareMathOperator{\diag}{diag}
\DeclareMathOperator{\End}{End}
\DeclareMathOperator{\id}{id}
\DeclareMathOperator{\im}{im}
\DeclareMathOperator{\Proj}{\mathbb P}
\DeclareMathOperator{\re}{Re}

\newcommand{\tS}{\mathscr S}

\newcommand{\hor}{{\mathscr H}}
\newcommand{\Hrel}{\mathrel{\sim_{\mkern-8mu\hor}^{}}}

\newcommand{\HKT}{\textsc{hkt}\xspace}
\newcommand{\KT}{\textsc{kt}\xspace}

\newcommand{\SKT}{\textsc{skt}\xspace}
\newcommand{\Hrelated}{\( \hor \)\bdash related\xspace}

\newcommand{\Horb}{H_{\textup{orb}}}
\newcommand{\LC}{\nabla^{\textup{LC}}}
\newcommand{\nB}{\nabla^{\textup{B}}}
\newcommand{\Ob}{\nabla^{\textup{Ob}}}
\newcommand{\TB}{T^{\textup{B}}}

\newcommand{\hook}{\mathop{\lrcorner}}
\newcommand{\ot}{\otimes}
\newcommand{\lift}[1]{\mathring{#1}}
\newcommand{\liftAM}{\lift A_M}
\newcommand{\liftTM}{\lift T_M}

\DeclarePairedDelimiter{\rcomp}{\llbracket}{\rrbracket}

\newcommand{\br}{\hspace{0pt}}
\newcommand{\bdash}{-\br} 
\newcommand{\eqbreak}{\\&\qquad}

\begin{document}
\begin{flushright}
  \begin{tiny}
    NORDITA-2008-61\\
    IMADA-PP-2008-16\\[5ex]
  \end{tiny}
\end{flushright}

\title{Twisting Hermitian and Hypercomplex Geometries}

\author{Andrew Swann}
\address[Swann]{Department of Mathematics and Computer Science\\
University of Southern Denmark\\
Campusvej 55\\
DK-5230 Odense M\\
Denmark} \email{swann@imada.sdu.dk}

\begin{abstract}
  A twist construction for manifolds with torus action is described
  generalising certain T-duality examples and constructions in hypercomplex
  geometry.  It is applied to complex, SKT, hypercomplex and HKT manifolds
  to construct compact simply-connected examples.  In particular, we find
  hypercomplex manifolds that admit no compatible HKT metric, and HKT
  manifolds whose Obata connection has holonomy contained in \(
  SL(n,\mathbb H) \).
\end{abstract}

\subjclass[2000]{Primary 53C55; Secondary 53C25, 53C29, 57S25, 32C37}

\keywords{Torus action, duality, Hermitian, hypercomplex, KT, SKT, HKT}

\maketitle

\begin{center}
  \begin{minipage}{0.6\linewidth}
    \begin{tiny}
      \tableofcontents
    \end{tiny}
  \end{minipage}
\end{center}

\newpage
\section{Introduction}
\label{sec:introduction}

The study of special metrics compatible with one or more complex structures
has a long history.  The most widely studied case is that of Kähler
metrics, where the Levi-Civita connection also preserves the complex
structure; here there is a rich plentiful source of examples, although
there are complex manifolds that admit no Kähler metric.  However, the
analogue of this situation for two anti-commuting complex structures,
namely hyperKähler geometry, is very restrictive and only a limited number
of compact examples are known.  Thus from the mathematical point of view,
weaker compatibility conditions are of interest.  On the other hand, models
from theoretical physics, particularly in the presence of supersymmetry,
lead to complex structures and metrics with potentially less restrictive
constraints \cite{Gates-HR:twisted,Strominger:superstrings}.  In general,
it is not hard to weaken the formal definitions, but then the question
remains whether these really lead new structures and how one might
construct examples.  In this paper, we present a general construction that
from a manifold with torus symmetry produces a new manifold, the twist, and
show how geometric data may be moved through this construction.  We then
specialise the construction to complex and hypercomplex geometry, using it
to construct examples that are compact and simply-connected.

In outline the twist construction we describe is as follows.  Consider a
manifold \( M \) with an action of an \( n \)-torus~\( \TM \).  Suppose \(
P\to M \) is a principal \( \Tn \)-bundle with connection.  If the \( \TM
\)-action lifts to \( P \) commuting with the principal action, then we may
construct the quotient space
\begin{equation*}
  W=P/\TM.
\end{equation*}
Furthermore, if the lifted \( \TM \)-action preserves the principal
connection, then tensors on \( M \) may be transferred to tensors on~\( W
\) by requiring their pull-backs to~\( P \) to agree on horizontal vectors.
In this way, an invariant geometric structure on~\( M \), such as complex
structure or a metric, determines a corresponding geometric structure on
the twist~\( W \).

Essentially this construction was considered mathematically by Joyce
\cite{Joyce:hypercomplex} in the context of hypercomplex manifolds using
instanton connections.  For the special case of circle actions and
instanton connections this was specialised to \HKT metrics for hypercomplex
structures by Grantcharov \& Poon \cite{Grantcharov-P:HKT}.  As
demonstrated in \cite{Swann:T}, specific non-compact examples in the \HKT
context include the T-duality construction evoked by Gibbons, Papadopoulos
\& Stelle \cite{Gibbons-PS:hkt-okt} based on the \( \sigma \)-model duality
of Buscher~\cite{Buscher:symmetry} (cf.\
\cite{Bergshoeff-HO:duality,Rocek-V:duality,Witten:various}).  The
construction we present is more general, applying in principle to any
geometric structure, but even in the hypercomplex or \HKT case we will see
that the instanton condition can be relaxed with advantage.  Furthermore,
in many situations the \( \Tn \)-twist construction is not equivalent to \(
n \) invocations of the \( S^1 \)-twist.  Versions of the \( S^1 \)-twist
have been announced and discussed in \cite{Swann:T,Swann:Bilbao} and
applied to almost quaternion-Hermitian manifolds
in~\cite{Cabrera-S:aqH-torsion}.

We start the paper by discussing the general framework for the twist
construction.  In particular, we study in detail the problem of lifting the
\( \TM \) action to the principal bundle~\( P \), \S\ref{sec:lift}, and
invoke the topological results of \cite{Lashof-MS:equivariant}.  The twist
construction itself is described in \S\ref{sec:twist}.  We show how tensors
and almost complex structures may be moved through the twist construction,
and study the effects on the exterior derivatives of forms and
integrability of complex structures.  In the latter case, we see the
instanton condition is not the most general requirement, in line with the
constructions of Goldstein \& Prokushkin \cite{Goldstein-P:SU3}.  As the
twist \( W \) is constructed as a quotient it is potentially singular.  For
our applications, we are only interested in smooth manifolds, but many of
the results extend without change to the orbifold case.

In \S\ref{sec:torsion}, we apply the twist construction in the context of
Hermitian geometry.  We concentrate mostly on the case of \SKT manifolds,
strong Kähler manifolds with torsion.  These are characterised by
\( \partial\bar\partial \omega_I = 0 \), where \( \omega_I \) is the Kähler
form.  Gauduchon~\cite{Gauduchon:canonical} showed that any compact
Hermitian surface admits such a metric, but in higher dimensions the
condition more complicated.  The \SKT structures on six-dimensional
nilmanifolds were classified by Fino, Parton \& Salamon \cite{Fino-PS:SKT},
however these are not simply-connected.  Grantcharov, Grantcharov \& Poon
\cite{Grantcharov-GP:CY-toric} used torus bundles to provide
six-dimensional examples.  Using the instanton case, we reproduce their
examples via the twist construction, extending them to higher dimensions,
and point out how a number of relatively explicit examples may be produced.
We also show that the non-instanton case gives rise to further compact
simply-connected \SKT manifolds.  This section closes with a discussion of
the behaviour of complex volume forms.  Note that other known compact
examples of \SKT manifolds include even-dimensional compact Lie groups
\cite{Spindel-STvP:complex} and certain instanton moduli spaces
\cite{Luebke-T:universal,Cavalcanti:metric-reduction}, and that Fino \&
Tomassini \cite{Fino-T:blow-up} have recently shown that the \SKT condition
is preserved by blow-up.

The final part of the paper, \S\ref{sec:HKT}, concerns hypercomplex and
\HKT geometry.  Hypercomplex manifolds carry two anti-commuting complex
structures; the \HKT condition on a compatible metric may be expressed via
a simple first order relation on the exterior derivatives of the
corresponding Kähler forms, equation~\eqref{eq:HKT}.  The \HKT condition
implies that there is a unique connection preserving the metric and the
complex structures, with torsion determined by a three-form.  This geometry
was originally introduced by Howe \& Papadopoulos
\cite{Howe-P:twistor-kaehler} in the physics literature and the paper of
Grantcharov \& Poon \cite{Grantcharov-P:HKT} gave a mathematical
description and several constructions.  We now know that \HKT structures
behave in many ways as a good quaternionic analogue of Kähler geometry.  In
particular, there is a potential theory~\cite{Banos-S:HKT}, a version of
Hodge theory~\cite{Verbitsky:Hodge} and some work towards Calabi conjecture
types of results~\cite{Aleskar-V:Calabi}.  However, the strongest versions
of these results require a reduction of the holonomy of the hypercomplex
structure: there is a unique torsion-free connection, the Obata connection,
that preserves the given complex structures, and the requirement is that
its holonomy should be contained the subgroup \( \SL(n,\bH) \), consisting
of the invertible \( n\times n \) quaternion matrices \( \GL(n,\bH)\subset
\GL(4n,\bR) \) that preserve a real volume form.

In \cite{Verbitsky:hyperholomorphic-kaehler}, Verbitsky used instanton
connections on vector bundles to construct compact \HKT manifolds, however
these have infinite fundamental group.  We show in
\S\ref{sec:instanton-HKT} that the instanton version of the twist
construction leads to many simply-connected examples.  We also demonstrate
that, at least locally, \HKT metrics may be produced from non-instanton
twists.  

Using instanton twists, we show in \S\ref{sec:SLnH}, that there are
non-trivial \HKT metrics on compact simply-connected manifolds such that
the Obata holonomy lies in~\( \SL(n,\bH) \).  Barberis, Dotti \& Verbitsky
\cite{Barberis-DV:canonical} recently provided similar examples on compact
nilmanifolds, but these have infinite fundamental group.

Finally, we construct via non-instanton twists examples of compact
simply-connected hypercomplex manifolds in all allowable dimensions that do
not admit a compatible \HKT metric.  Examples with infinite fundamental
group were previously constructed on nilmanifolds in dimension~\( 8 \),
Fino \& Grantcharov \cite{Fino-G:SU-Sp}.

\begin{acknowledgements}
  I thank Martin Svensson, Francisco Martín Cabrera, Andrew Dancer and Gil
  Cavalcanti for useful conversations at various stages of this work that
  is partially supported by a grant from the MEC (Spain), project
  MTM2004-2644.  It is also a pleasure to thank NORDITA and the organisers
  of the programme on `Geometrical Aspects of String Theory' for
  support during the completion of this paper.
\end{acknowledgements}

\section{Lifting Abelian actions}
\label{sec:lift}

Let \( \An \) be a connected Abelian Lie group of dimension~\( n \).
Suppose that \( \pi\colon P\to M \) is a principal \( \An \)-bundle with
structural group~\( \AP \).  Write \( \aP \) for the Lie algebra of~\( \AP
\) and let \( \rho\colon \aP \to \XM(P) \), \( Y\mapsto \rho_Y \), be the
vector fields generated by the principal action.

We now assume there is an action of \( \AM\cong \An \) on~\( M \) and write
\( \xi\colon \aM \to \XM(M) \) for the infinitesimal action.  If \( \theta
\in \Omega^1(P,\aP) \) is a connection one-form on~\( P \) with curvature
\( F \in \Omega^2(M,\aP) \), \( \pi^*F=d\theta \), we wish to determine
conditions so that the \( \AM \)-action is covered by an Abelian Lie group
action on~\( P \) preserving~\( \theta \) and commuting with~\( \AP \).  We
will use \( \tilde X \) to denote the horizontal lift of \( X\in TM \)
to~\( \hor = \ker \theta \subset TP \).

First consider the problem of lifting to an \( \bR^n=\widetilde{\An}
\)-action.  Such a lift is given by a map \( \lift\xi\colon \aM \to \XM(P)
\).

\begin{proposition}
  \label{prop:Rn-lift}
  The \( \AM \)-action induced by \( \xi \) lifts to an \( \bR^n \)-action
  preserving the connection form~\( \theta \) if and only if
  \begin{enumerate}
  \item\label{item:lift1} \( L_\xi F=0 \),
  \item\label{item:lift2} \( [\xi \hook F]=0\in H^1(M,\aP\ot\aM^*) \) and
  \item\label{item:lift3} \( \xi^*F = 0 \).
  \end{enumerate}
  Moreover, if \( \AM \) is compact, then condition~\itref{item:lift3} is
  redundant.
\end{proposition}

\begin{proof}
  Write
  \begin{equation}
    \label{eq:lift}
    \lift\xi = \tilde\xi + \lift a\rho
  \end{equation}
  for some \( \lift a\in \Omega^0(P,\aP\ot\aM^*) \), where \( \tilde\xi \)
  is the horizontal lift and \( \rho \)~is regarded as an element of~\(
  \XM(P,\aP^*) \).  The condition that \( \lift\xi \) preserves \( \theta
  \) gives
  \begin{equation*}
    0 = L_{\lift\xi}\theta
    = \lift\xi \hook d\theta + d(\lift\xi \hook \theta)
    = \lift\xi \hook \pi^*F + d(\lift a\rho \hook \theta)
    = \pi^*(\xi\hook F) + d\lift a.
  \end{equation*}
  Thus \( d\lift a = -\pi^*(\xi\hook F) \).  The first consequence of this
  is that \( \pi^*(L_\xi F) = \pi^* (d\xi\hook F) = - d^2\lift a = 0 \),
  giving condition~\itref{item:lift1}.  We also find differentiating in the
  vertical directions that \( \rho\lift a =
  \rho \hook d\lift a = - \rho \hook \pi^*d(\xi \hook F) = 0 \), so \(
  \lift a \) is constant on fibres and is the pull-back of a function \(
  a\in\Omega^0(M,\aP\ot\aM^*) \):
  \begin{equation*}
    \lift\xi\hook\theta = \lift a = \pi^*a.
  \end{equation*}
  We now have
  \begin{equation}
    \label{eq:a}
    da = - \xi \hook F
  \end{equation}
  so the class \( [\xi\hook F] \) is zero in \( H^1(M,\aP\ot\aM^*) \),
  which is condition~\itref{item:lift2}.

  It remains to show that \( \lift\xi\colon \aM \to \XM(P) \) is a Lie
  algebra homomorphism.  As we are working with Abelian groups, this is the
  same as \( [\lift\xi_X,\lift\xi_Y]=0 \) for all \( X,Y\in\aM \).  We have
  \begin{equation}
    \label{eq:lift-br}
    [\lift\xi_X,\lift\xi_Y] = [\tilde\xi_X,\tilde\xi_Y] + [\lift
    a_X\rho,\tilde\xi_Y] + [\tilde\xi_X,\lift a_Y\rho] + [\lift
    a_X\rho,\lift a_Y\rho]. 
  \end{equation}
  The last term vanishes since \( \aP \) is Abelian and \( \rho\lift a=0
  \).  For one of the middle terms, we have
  \begin{equation*}
    \pi_*[\lift a_X\rho,\tilde\xi_Y] = [0,\xi_Y] = 0,
  \end{equation*}
  so the horizontal part is zero, and the vertical part is determined by
  \begin{equation*}
    \begin{split}
      \theta([\lift a_X\rho,\tilde\xi_Y])
      & = -d\theta(\lift a_X\rho,\tilde\xi_Y) - \tilde\xi_Y(\lift a_X) \\
      & = -(\pi^*F)(\lift a_X\rho,\tilde\xi_Y) - \pi^*(\xi_Ya_X) \\
      & = \pi^*(F(\xi_X,\xi_Y)).
    \end{split}
  \end{equation*}
  As \( \aM \) is Abelian, the first term in \eqref{eq:lift-br} has
  \begin{equation*}
    \pi_*([\tilde\xi_X,\tilde\xi_Y]) = [\xi_X,\xi_Y] =
    0\quad\text{and}\quad 
    \theta([\tilde\xi_X,\tilde\xi_Y]) = -\pi^*(F(\xi_X,\xi_Y)).
  \end{equation*}
  Putting these together we get
  \begin{equation*}
    \theta([\lift\xi_X,\lift\xi_Y]) = \pi^*(F(\xi_X,\xi_Y)),
  \end{equation*}
  whilst the horizontal part is zero.  This gives
  condition~\itref{item:lift3}.  Now note that \( F(\xi_X,\xi_Y)=
  -da_X(\xi_Y) = - L_{\xi_Y}a_X \), so \( d(F(\xi_X,\xi_Y)) = -
  L_{\xi_Y}da_X = L_{\xi_Y}(\xi_X\hook F) = 0 \).  This shows that \(
  F(\xi_X,\xi_Y) \) is constant.  When \( \AM \)~is compact, each component
  of \( \xi_Ya_X \) has a zero on each \( \AM \)-orbit, e.g., at a maximum
  of the component of \( a_X \), so the constant \( F(\xi_X,\xi_Y) \) is
  zero.
\end{proof}

Note that the lift is not unique, depending instead on a choice of \(
a\in\Omega^0(M,\aP\ot\aM^*) \) in equation~\eqref{eq:a}.  In particular, if
\( M \) is compact, we can add a constant element of \(\aP\ot\aM^* \) to~\(
a \) to ensure that \( a \)~is invertible.  

Suppose \( F \) is a closed \( 2 \)-form with values in \( \mathbb
R^n\cong\aP \).

\begin{definition}
  \label{def:F-Ham}
  We say that an \( \AM \)-action is \emph{\( F \)-Hamiltonian} if it
  satisfies conditions \itref{item:lift1} and~\itref{item:lift2} of
  Proposition~\ref{prop:Rn-lift}.
\end{definition}

\subsection{Topological considerations for torus actions}
\label{sec:topology}

In constructions with \( \An=\Tn \) an \( n \)-torus, we would like have as
starting data a manifold \( M \) with \( \TM \)-action and a \( 2 \)-form
\( F\in \Omega^2(M,\tn) \) that is \( \TM \)-invariant.  We then wish to
construct a principal \( \Tn \)-bundle~\( P \) with a connection \( \theta
\) whose curvature is \( F \) in such a way that \( \TM \) lifts to a \(
\Tn \)-action on \( P \) preserving \( \theta \) and commuting with the
principal action.  

First, ignoring the \( \TM \)-action we need \( F \)~to be a closed form
with integral periods, we write this as \( F\in\Omega^2_\bZ(M,\tn) \).  We
then have that \( [F] \in H^2(M,\bZ^n)\ot \bR\subset H^2(M,\tn) \) and by
Chern-Weil theory we can find a principal \( \Tn \)-bundle \( P \) with \(
c_1(P)\ot\bR = [F] \) and connection \( \theta_0 \) on \( P \) with \(
d\theta_0 = \pi^*F \).  Let us write \( \Tn = \tn/\Ln \) and \( \TM =
\tM/\LM \), so \( \Ln\cong \bZ^n \cong \LM \).

The question of equivariantly lifting the \( \TM \)-action to~\( P \) is
mostly clearly addressed in \cite{Lashof-MS:equivariant}.  The main tool
used is the spectral sequence of the fibration
\begin{equation*}
  M \longrightarrow M_G \longrightarrow BG,
\end{equation*}
where \( G=\TM \), \( EG=(S^\infty)^n=(\colim S^k)^n\to BG=B\TM \cong
\CP(\infty)^n \) is the classifying space and \( M_G = EG \times_G M \).
It is shown that the \( \TM \)-action lifts to \( P \) if and only if \(
c_1(P)\in H^2(M,\Ln) \) is \( \TM \)-invariant and lies in
\begin{equation*}
  E_\infty^{0,2} = \ker(d_3\colon E_3^{0,2}\to E_3^{3,0}) \subset
  \ker(d_2\colon E_2^{0,2}\to E_2^{2,1}).
\end{equation*}
Since \( \TM \) is connected it acts trivially on \( H^q(M,\Ln) \), so the
invariance of \( c_1(P) \) is automatic, and the locally constant
presheaf~\( H^q(M,\Ln) \) over the simply-connected space~\( BG \) is
constant.  We thus have
\begin{equation*}
  \begin{split}
    E_2^{p,q}
    &= H^p(B\TM,H^q(M,\Ln))\\
    &=
    \begin{cases}
      H^q(M,\Ln)\ot_\bZ S^k\LM^*,&\text{for \( p=2k\geqslant0 \)
      even,}\\
      0,&\text{otherwise.}
    \end{cases}
  \end{split}
\end{equation*}
This gives immediately that \( d_3 = 0 \).  Since \( E_3^{p,q} =
\ker(d_2\colon E_2^{p,q}\to E_2^{p+2,q-1}) / \im(d_2\colon E_2^{p-2,q+1}\to
E_2^{p,q}) \), we have
\begin{equation*}
  \begin{split}
    E_\infty^{0,2}
    &= E_3^{0,2} = \ker(d_2\colon E_2^{0,2}\to E_2^{2,1})\\
    &= \ker (d_2\colon H^2(M,\Ln) \to H^1(M,\Ln)\ot_\bZ\LM^*).
  \end{split}
\end{equation*}
Tensoring with \( \bR \), the map \( d_2 \) is given by
\begin{equation*}
  (d_2^{\bR})[F] = [\xi\hook F] \in H^1(M,\tn\ot\tM^*) 
\end{equation*}
when \( F \)~is \( \TM \)-invariant.

If the \( \TM \)-action is \( F \)-Hamiltonian
(Definition~\ref{def:F-Ham}), then we have \( [\xi\hook F]=0 \) and the
class \( d_2(c_1(P)) \) is torsion in \( H^1(M,\Ln)\ot_\bZ\LM^* \).
However, \( H^1(M,\bZ) \) is isomorphic to the torsion-free part of \(
H_1(M,\bZ) \), so \( d_2(c_1(P)) = 0 \).  Thus the \( \TM \)-action lifts
to~\( P \) covering the action on~\( M \).  The lifted \( \TM \)-action
will not necessarily preserve \( \theta_0 \), but averaging \( \theta_0 \)
over the lifted \( \TM \) gives an invariant connection form \( \theta \)
still with curvature~\( F \).

In summary, we have:

\begin{proposition}
  Suppose \( M \) admits an \( F \)-Hamiltonian \( \TM \)-action for some
  closed \( 2 \)-form with integral periods \( F\in\Omega^2_{\mathbb
  Z}(M,\tn) \).  Then there is a principal \( \Tn \)-bundle \( P\to M \)
  which admits
  \begin{enumerate}
  \item a \( \liftTM \)-action of the \( n \)-torus~\( \TM \) commuting
    with the principal action, covering the \( \TM \)-action on~\( M \),
    and
  \item a \( \liftTM \)-invariant \( \Tn \)-connection \( \theta \) on~\( P
    \) with curvature~\( F \).
  \end{enumerate}
  In fact, such a lift exists for any \( \Tn \)-bundle \( P \) with \(
  c_1(P)\otimes \bR = [F] \).\qed
\end{proposition}

The lifts above are not unique: when lifting to an \( \bR^n \)-action,
there is an ambiguity in \( a\in\Omega^0(M,\tP\ot\tM^*) \) as one may add
any constant element of \( \lambda\in\tP\ot\tM^* \) to~\( a \).  Given a
torus lift, the integrality condition \( \lambda \in \LP\ot\LM^*\ot\bZ \)
leads to a new such torus action on~\( P \).  If the \( \TM \)-action on~\(
M \) is free, then each of these lifts \( \liftTM \) is free on~\( P \).
Allowing \( \lambda \in \LP\ot\LM^*\ot\bQ \) leads to torus actions on~\( P
\) whose infinitesimal generators map to those on~\( M \) under \( \pi_* \)
and which cover the action of a finite cover of~\( \TM \).

\section{Twist construction}
\label{sec:twist}

Suppose \( M \) is a manifold with an effective \( F \)-Hamiltonian \( \AM
\)-action where \( F\in\Omega^2_\bZ(\an) \).  Let \( (P,\theta) \) be a
principal \( \An \)-bundle with curvature~\( F \) and an \( \liftAM
\)-action preserving \( \theta \) and covering the \( \AM \)-action
infinitesimally.  Here \( \liftAM \) is some connected Abelian group
covering~\( \AM \).  Assume that \( \liftAM \) acts properly on~\( P \);
and that \( \liftAM \) is transverse to
\begin{equation*}
  \hor = \ker\theta.
\end{equation*}
Then \( \liftAM \) has discrete stabilisers and \( P/\liftAM \) has the
same dimension as~\( M \).  This transversality is the same as requiring \(
a\in\Omega^0(M,\aP\ot\aM^*) \) in~\eqref{eq:a} to be invertible.  If \( M
\) is compact and \( \An\cong\Tn \) is a torus, then the discussion of the
previous section shows that there is always a proper lift and that we may
add a constant rational element of \( \aP\ot\aM^* \) to~\( a \) to ensure
that \( a^{-1}\in\Omega^0(M,\aM\ot\aP^*) \) exists.

\begin{definition}
  A \emph{twist} of \( M \) with respect to \( \AM \), \( F \) and
  invertible~\( a \), is the quotient space
  \begin{equation*}
    W = P / \liftAM.
  \end{equation*}
  We say that \( W \) is a \emph{smooth twist} if \( W \) is a manifold.
\end{definition}

For torus actions, a twist~\( W \) will at worst be an orbifold under the
assumptions above.  We are interested in constructing smooth manifolds and
will therefore only discuss geometric structures in the case of smooth
twists, however many of our results will carry over to the orbifold case
without modification.  Note that the following example shows that smooth
twists may not exist.

\begin{example}
  Suppose \( M=\CP(n) \) with circle action
  \begin{equation*}
    [z_0,\dots,z_n]\mapsto[e^{2\pi i\lambda_0\theta}z_0,\dots,e^{2\pi
    i\lambda_n\theta}z_n]. 
  \end{equation*}
  The canonical circle bundle over \( \CP(n) \) is \( S^{2n+1}\subset
  \bC^{n+1} \) with principal action \( \mathbf z \mapsto e^{2\pi i\phi}\mathbf
  z \).  The lifted action
  \begin{equation*}
    (z_0,\dots,z_n)\mapsto(e^{2\pi i(\lambda_0+k)t}z_0,\dots,e^{2\pi
    i(\lambda_n+k)t}z_n) 
  \end{equation*}
  is free only if \( k=-\tfrac12(\lambda_i+\lambda_j) \) for all \(
  \lambda_i\ne\lambda_j \), which in general is impossible.
\end{example}

In the case that the lifting function \( a \) is non-constant, we can not
consistently identify \( \tM \) and \( \tP \), and so twisting with \( \TM
\) can not be reduced to repeated twists by circle subgroups via circle
bundles.

\subsection{Geometric structures}
\label{sec:geometry}

Returning to \( W \) a smooth twist of a general \( M \), we have
projection maps
\begin{equation*}
  \begin{CD}
    M @<\pi<< P @>\pi_W^{}>> W.
  \end{CD}
\end{equation*}
Our assumptions imply that both maps are transverse to the distribution~\(
\hor \).  We use this to relate objects on \( M \) and~\( W \).  Invariant
vector fields may be transferred simply by lifting horizontally and
projecting.

\begin{definition}
  Two \( (p,0) \)-tensors \( \alpha \) on \( M \) and \( \alpha_W \) on \(
  W \) are said to be \emph{\Hrelated{}}, written
  \begin{equation*}
    \alpha \Hrel \alpha_W,
  \end{equation*}
  if their pull-backs to \( P \) agree on \( \hor \), i.e. \( \pi^*\alpha =
  \pi_W^*\alpha_W^{} \) on~\( \hor \).
\end{definition}

\begin{lemma}
  \label{lem:H-p}
  Each invariant \( p \)-form \( \alpha\in\Omega^p(M)^{\AM} \) is \Hrelated
  to a unique \( p \)-form \( \alpha_W \in \Omega^p(W) \) given by
  \begin{equation*}
    \pi_W^*\alpha_W^{} = \pi^*\alpha - \theta \wedge
    \pi^*(a^{-1}\xi\hook\alpha).  
  \end{equation*}
\end{lemma}

Note that \( \AM \)-invariance of~\( \alpha \) is a necessary condition.

\begin{proof}
  The form \( \pi_W^*\alpha_W^{} \) may be decomposed with respect to \(
  \theta \) and \( \hor \) as an element of \( \Omega^p(\hor) + \theta
  \wedge \Omega^{p-1}(\hor,\aP^*) \).  By definition the first component is
  \( \pi^*\alpha \); we write the second as \( \theta \wedge \beta \).
  Using \eqref{eq:lift}, we now compute
  \begin{equation*}
    \begin{split}
      0 &= \lift\xi\hook\pi_W^*\alpha_W^{} = \lift\xi\hook\pi^*\alpha +
      \lift\xi\hook(\theta \wedge \beta) \\
      &= \tilde\xi\hook\pi^*\alpha + \lift a\theta(\rho)\beta -
      \theta\wedge\tilde\xi\hook\beta \\
      &= \pi^*(\xi\hook\alpha) + (\pi^*a)\beta - \theta \wedge
      (\xi\hook\beta),
    \end{split}
  \end{equation*}
  since \( \theta\circ\rho \) is the identity on~\( \aP \).  Considering
  horizontal vectors, we have \( \beta = -\pi^*(a^{-1}\xi\hook\alpha) \)
  and the result follows.
\end{proof}

\begin{corollary}
  \label{cor:g}
  If \( g \)~is an invariant metric on~\( M \) then the unique metric~\(
  g_W \) on~\( W \) \Hrelated to~\( g \) is given by
  \begin{equation*}
    \pi_W^*g_W^{} = \pi^*g - 2\theta \vee \pi^*(a^{-1}\xi^\flat) +
    \theta^2\pi^*((a^{-1})^2g(\xi\ot\xi)).
  \end{equation*}
  \qed
\end{corollary}

\begin{corollary}
  \label{cor:d}
  Suppose \( \alpha \Hrel \alpha_W \).  Then
  \begin{equation*}
    d\alpha_W \Hrel d\alpha -  a^{-1}F\wedge\xi\hook\alpha.
  \end{equation*}
\end{corollary}

\begin{proof}
  Computing directly, we have
  \begin{equation*}
    \begin{split}
      \pi_W^*d\alpha_W^{} &= d\pi_W^*\alpha_W^{} = d(\pi^*\alpha - \theta
      \wedge \pi^*(a^{-1}\xi\hook\alpha))\\
      &= \pi^*(d\alpha - F\wedge a^{-1}\xi\hook\alpha) +
      \theta\wedge\pi^*(d(a^{-1}\xi\hook\alpha))
    \end{split}
  \end{equation*}
  which agrees horizontally with the claimed result.  Vertically we have \(
  \theta \) wedge the pull-back under \( \pi \) of
  \begin{equation*}
    \begin{split}
      -a^{-1}{da}\,a^{-1}\wedge \xi\hook\alpha + a^{-1} d(\xi\hook\alpha)
      &= a^{-1}\xi\hook F\,a^{-1}\wedge \xi\hook\alpha - a^{-1}
      \xi\hook d\alpha)\\
      &= - (a^{-1}\xi\hook(d\alpha - a^{-1}F\wedge \xi\hook\alpha) ),
    \end{split}
  \end{equation*}
  since \( L_\xi\alpha=0 \).  The result now follows.
\end{proof}

Thus the invariant part of the exterior algebra of~\( W \) may be regarded
as the invariant exterior algebra of~\( M \) with the twisted differential
\( d - a^{-1}F\wedge\xi\hook \).

\subsection{Duality}
\label{sec:duality}

Let us now show that \( W \) is dual to \( M \) in the sense that \( M \)
may also be obtained from \( W \) via a twist.  The distribution \( \hor \)
on~\( P \) is transverse to the action of \( \liftAM \).  If \( \liftAM \)
acts freely we have a principal bundle \( \liftAM\to P\to W \).  Its
connection form corresponding to~\( \hor \) is
\begin{equation*}
  \theta_W = \pi^*(a^{-1})\theta,
\end{equation*}
as may be seen by writing \( \theta_W = f\theta \), for some \(
f\in\Omega^0(P,\aM\ot\aP^*) \), and enforcing the condition \(
\theta_W(\lift\xi) = \id \) on~\( \aM \).  This has curvature
\begin{equation*}
  \pi_W^*F_W^{} = \pi^*(a^{-1}F) - \pi^*(a^{-1}da\, a^{-1})\wedge\theta,
\end{equation*}
which is simply the two-form \( F_W \) that is \( \hor \)-related to~\( F
\).  Since \( \AP \) commutes with \( \liftAM \), it descends to an action
of an Abelian group \( \AW \) on~\( W \) preserving \( F_W \).  Write \(
\zeta\colon \aW \to \XM(W) \) for the infinitesimal action of~\( \AW \), so
\( \zeta = \pi_W\circ \rho \).  This action is \( F_W \)-Hamiltonian with
\( \zeta\hook F_W = -d(a^{-1}) \).  The original manifold \( M \) is
obtained by twisting \( W \) with respect to \( \AW \) and \( F_W \) using
\( a^{-1} \).

\subsection{Lie brackets and complex structures}
\label{sec:Lie}

Tangent vectors \( X \) on~\( M \) and \( X_W \) on~\( W \) are said to be
\Hrelated if their horizontal lifts to \( \hor \) agree.  Writing \(
\widehat\cdot \) for the horizontal lift from~\( W \), this says,
\begin{equation*}
  \widehat{X_W} = \widetilde X.
\end{equation*}

\begin{lemma}
  \label{lem:Lie}
  Lie brackets between \Hrelated vector fields are related by
  \begin{equation*}
    [X_W,Y_W] \Hrel [X,Y] - \xi a^{-1} F(X,Y).
  \end{equation*}
\end{lemma}

\begin{proof}
  Lifting \( [X,Y] \) horizontally gives the horizontal part of \(
  [\widetilde X,\widetilde Y] \).  The vertical part of this last Lie
  bracket is \( \rho \pi^* F(X,Y) \), so
  \begin{equation*}
    \widetilde{[X,Y]} = [\widetilde X,\widetilde Y] - \rho\pi^* F(X,Y).
  \end{equation*}
  Similarly,
  \begin{equation*}
    \begin{split}
      \widehat{[X_W,Y_W]} &= [\widehat X_W,\widehat Y_W] - \lift\xi
      \pi_W^*F_W(X_W,Y_W) \\
      &= [\widehat X_W,\widehat Y_W] - (\tilde\xi+a\rho)\pi^*(a^{-1}F +
      a^{-1}da\, a^{-1}\wedge\theta)(\widehat{X_W},\widehat{Y_W})\\
      &= [\widehat X_W,\widehat Y_W] - (\pi^*(\xi a^{-1}F) +
      \rho\pi^*F)(\widehat{X_W},\widehat{Y_W}),
    \end{split}
  \end{equation*}
  from which the result follows.
\end{proof}

Almost complex structures may now be \Hrelated in a similar way giving
\begin{equation}
  \label{eq:I}
  \widehat{I_WA}= \widetilde{I\pi_*\hat A} .
\end{equation}

For a \( p \)-form \( \alpha \) and an index \( k \) we write
\begin{equation*}
  I_{(k)}\alpha(X_1,\dots,X_p) = -\alpha(X_1,\dots,IX_k,\dots,X_p) 
\end{equation*}
and \( I_{(ab\dots c)}=I_{(a)}I_{(b)}\dots I_{(c)} \).  This convention
ensures that \( I_{(k)}J_{(k)}\alpha = (IJ)_{(k)}\alpha \), but is the
opposite of the usual convention in complex geometry.  In particular, \(
\Lambda^{1,0} \) becomes the \( (-i) \)-eigenspace for \( I_{(1)} \).

\begin{proposition}
  \label{prop:N}
  The Nijenhuis tensors of \Hrelated almost complex structures \( I \) and
  \( I_W \) are related by
  \begin{equation}
    \label{eq:N}
    N_{I_W} \Hrel N_I  + (1-\cL_I)\cF,
  \end{equation}
  where \( \cF = \xi a^{-1} F \in \Gamma(TM\ot \Lambda^2T^*M) \) and
  \( \cL_I = I_{(12)}+I_{(13)}+I_{(23)} \).
\end{proposition}

\begin{proof}
  This follows directly from the definition \( N_I(X,Y)=
  [IX,IY]-I[IX,Y]-I[X,IY]-[X,Y] \) and Lemma~\ref{lem:Lie}.
\end{proof}

Note that \( \cL_I \) acts on \( TM\ot \Lambda^2T^*M \) with
eigenvalues \( -3 \) and \( +1 \).  The \( (-3) \)-eigenspace is \(
\rcomp{T^{1,0}\ot \Lambda^{0,2}} \), where \( \rcomp V \otimes \bC = V +
\overline V \).  We thus see that if \( F \) is of type~\( (1,1) \), then
twisting preserves integrability.  However, it will be important for us
that other choices of \( F \) can also give integrable complex structures.

To understand the integrability better, fix a point \( x\in M \), let \(
\sA=\im\xi\subset T_xM \) and put \( \sA_I= \sA\cap I\sA \).  Define \(
s=\dim_\bC\sA_I \) and \( r=\dim_\bR\sA \).  For a basis \( e_1,\dots,e_n
\) of~\( \aM \) write \( X_i=\xi_x(e_i) \).  Using \( \sA_I\leqslant
\sA\leqslant T_xM \), we may choose this basis so that
\begin{enumerate}
\item \( X_1,\dots,X_{2s} \) is a complex basis of \( \sA_I \), with \(
  IX_{2j-1}=X_{2j} \), for \( j=1,\dots,s \),
\item \( X_1,\dots,X_r \) is a basis of \( \sA \), and
\item \( X_{r+1},\dots,X_n \) are zero.  
\end{enumerate}
We then write \( a^{-1}F\in \Omega^2(M,\aM) \) as
\begin{equation*}
  a^{-1}F = \sum_{k=1}^n F_k\ot e_k
\end{equation*}
with \( F_k\in\Omega^2(M) \), so \( \cF = \sum_{k=1}^r X_k\ot F_k \).

\begin{lemma}
  \label{lem:integrable}
  If \( (M,I) \) is complex, then the induced almost complex structure \(
  I_W \) on a twist \( W \) is integrable if and only if
  \begin{enumerate}
  \item \( (F_{2j-1}+iF_{2j})^{0,2} =0 \), for \( j=1,\dots,s \), and
  \item \( F_k \in \Lambda^{1,1} \), for \( k=2s+1,\dots,r \).
  \end{enumerate}
\end{lemma}

\begin{proof}
  By \eqref{eq:N}, the integrability condition for \( I_W \) is \(
  (1-\cL_I)\cF = 0 \).  Note that
  \begin{equation*}
    \begin{split}
      (1-\cL_I)(X_k\ot F_k) 
      &= X_k\ot(1-I)F_k - IX_k\ot I_{(1)}(1-I)F_k\\
      &= 2X_k\ot (F_k^{2,0}+F_k^{0,2}) + 2IX_k\ot (iF_k^{2,0}-iF_k^{0,2}).
    \end{split}
  \end{equation*}
  For \( 2s<k\leqslant r \), \( X_k \) and \( IX_k \) are linearly
  independent of \( X_j \) for all \( j\ne k \), so \(
  F_k^{2,0}=0=F_k^{0,2} \).  For \( j\leqslant s \), we should consider the
  components \( 2j-1 \) and \( 2j \) together.  Let \( X=X_{2j-1} \), \(
  F_{(1)}=F_{2j-1} \), \( F_{(2)} = F_{2j} \) and \( F_{(c)} =
  F_{(1)}+iF_{(2)} \).  Now, since \( F_{(1)} \) and \( F_{(2)} \) are real,
  we get
  \begin{equation*}
    \begin{split}
      (1-\cL_I)&(X_{2j-1}\ot F_{2j-1} + X_{2j}\ot F_{2j}) \\
      &= 2X\ot (F_{(1)}^{2,0} + F_{(1)}^{0,2} - iF_{(2)}^{2,0} +
      iF_{(2)}^{0,2}) \eqbreak + 2IX\ot (F_{(2)}^{2,0} + F_{(2)}^{0,2} +
      iF_{(1)}^{2,0} - iF_{(1)}^{0,2}) \\
      &= 4\re((X-iIX)\ot F_{(c)}^{0,2})
    \end{split}
  \end{equation*}
  and the result follows.
\end{proof}

\section{Hermitian and SKT structures}
\label{sec:torsion}

Suppose \( (g,I) \) is a Hermitian structure on~\( M \), meaning that \( I
\) is integrable and \( g(IX,IY)=g(X,Y) \) for all \( X,Y\in TM \).  Then
there is a unique connection~\( \nB \) with skew-symmetric torsion
preserving \( I \) and~\( g \) \cite{Gauduchon:Hermitian-Dirac}.  This is
known as the Bismut connection, due to its appearance
in~\cite{Bismut:loc-nK}, and is given by
\begin{equation}
  \label{eq:Bismut}
  \nB = \LC + \tfrac12\TB,\qquad c = (\TB)^\flat = -Id\omega_I,
\end{equation}
where \( \omega_I(X,Y)=g(IX,Y) \).

\begin{definition}
  The triple \( (g,I,\nB) \), or equivalently \( (g,I,c) \), is known as a
  \emph{\KT structure}.  It is \emph{strong} or \emph{\SKT{}} if \( dc=0 \).
\end{definition}

Let us consider how the torsion form~\( c \) of a \KT manifold changes
under a twist.  We use the notation of~\S\ref{sec:twist} so that \(
\AM\cong\An \) is a connected Abelian group acting on~\( M \) in such a way
that there is a smooth twist~\( W \) given via curvature \(
F\in\Omega^2_\bZ(M,\aP) \), \( \AP\cong\An \), and lifting function~\(
a\in\Omega^0(M,\aP\ot\aM^*) \).

\begin{proposition}
  \label{prop:KT}
  Suppose the action of \( \AM \) preserves a Hermitian structure \( (g,I)
  \) on \( M \) with torsion form~\( c \).  If \( (F,a) \) are such that
  the twist has \( I_W \) integrable, then \( W \) is a \KT manifold with
  torsion given by
  \begin{equation}
    \label{eq:c}
    c_W \Hrel c - a^{-1}IF\wedge \xi^\flat.
  \end{equation}
\end{proposition}

\begin{proof}
  By construction \( I_W \) is compatible with \( g_W \).  We thus have \(
  \omega_I^W \Hrel \omega_I^{} \) and that \( c \) is invariant under \(
  \AM \).  Now the result follows from
  \begin{equation*}
    c_W \Hrel -I(d-a^{-1}F\wedge\xi\hook)\omega_I,
  \end{equation*}
  by Corollary~\ref{cor:d}.
\end{proof}

Many examples of this construction may be given by considering any free \(
\Tn \)-action on a Hermitian manifold and choosing a Hamiltonian two-form
\( F\in\Omega^2_\bZ(M,\bR^n) \).  As we will see below, this often produces
Hermitian structures that are not of Kähler type.

It follows from \eqref{eq:c} that the exterior derivative of the torsion
satisfies
\begin{equation}
  \label{eq:dc}
  \begin{split}
    dc_W &\Hrel dc - a^{-1}\bigl(F\wedge(\xi\hook c) + d(IF)\wedge
    \xi^\flat + IF \wedge d\xi^\flat\bigr) \eqbreak +
    (a^{-1})^2\bigl(g(\xi,\xi)F\wedge IF + F\wedge(\xi\hook IF)\wedge
    \xi^\flat \\ &\hspace{8em} -(\xi\hook F)\wedge IF\wedge
    \xi^\flat\bigr).
  \end{split}
\end{equation}
Note that if \( F \) is of type \( (1,1) \), the \emph{instanton case},
then this simplifies to
\begin{equation}
  \label{eq:dc11}
  dc_W \Hrel dc - a^{-1}F\wedge\bigl((\xi\hook c)  + d\xi^\flat -
  g(\xi,\xi)a^{-1}F \bigr).
\end{equation}

Using these last two expressions it is now possible to give a number of
examples of compact simply-connected \SKT manifolds.

\subsection{Instanton twists from tori}
\label{sec:SKT}

First note that any Kähler manifold is \SKT since it has \( c=0 \).  Let \(
N \) be an \SKT manifold and consider the product \( M = N\times \Tn[2] \),
where \( \Tn[2] \) is a \( 2 \)-torus with an invariant, so flat, Kähler
metric.  Then \( M \) is \SKT with torsion supported on~\( N \).  Let \(
\xi \) be the torus action on the \( \Tn[2] \)-factor.  We have \(
\xi\hook c=0 \) and \( d\xi^\flat=0 \).  If \( F\in\Omega^2_\bZ(N,\aP) \) with
\( \aP \cong \bR^2 \), then \( \xi\hook F=0 \) and \( a \)~is a constant
isomorphism \( \aM\to\aP \).  Choosing bases we may write \( \xi =
(X_1,X_2) \), \( a^{-1}F = (F_1,F_2) \).

Suppose now that \( F \) is of type~\( (1,1) \).  Then the twist \( W \) is
Hermitian and the \SKT condition becomes
\begin{equation*}
  \sum_{i,j=1}^2 g(X_i,X_j)F_i\wedge F_j=0.
\end{equation*}
Taking \( \Tn[2] = \bC/(\bZ + i\bZ) \) with the standard metric, we may
choose \( X_i \) orthonormal and \( a \) the identity matrix, the strong
condition reduces to \( F_1^2+F_2^2=0 \) with \( F_i\in\Omega^{1,1}_\bZ(N)
\).

\begin{example}
  Let \( N = \Sigma_1\times \Sigma_2 \) a product of Riemann surfaces and
  take \( [F_i] \) to be the fundamental class on the \( i \)th factor.
  The twist \( W \) is then a product of non-trivial circle bundles over
  the Riemann surfaces.  When each \( \Sigma_i \) has genus~\( 0 \), the
  twist \( W \) is topologically the product \( S^3\times S^3 \) and we
  obtain the Calabi-Eckmann complex structures~\cite{Calabi-E:complex}.
  The torsion form \( c_W \) is a sum of volume forms on the two factors.
  In general, the constructed \SKT structure only has \( \Tn[2]
  \)-symmetry, and even that may be destroyed by adding \(
  i\partial\bar\partial f \) to the Kähler form \( \omega_I^W \) for \(
  f\in C^\infty(W) \) with suitably small \( C^2 \)-norm.
\end{example}

To construct further examples of this type, note that Eells \& Lemaire
\cite{Eells-Lemaire:CBMS} showed that for an almost Kähler manifold~\( X \),
any non-constant holomorphic map \( f\colon X\to\CP(1) \) has \(
[f^*\omega_{\CP(1)}]\ne 0 \) in \( H^2(X,\bZ) \).  Thus if \( N \) is
Kähler and admits a non-constant holomorphic map to \( \CP(1) \), we get a
two-form
\begin{equation*}
  \Phi(f)=2\pi f^*\omega_{\CP(1)}\in\Omega_\bZ^{1,1}(N) 
\end{equation*}
with \( \Phi(f)^2=0 \).  Furthermore such a map \( f \) will exist whenever
\( N \) admits a non-constant holomorphic map to some compact Riemann
surface.

\begin{example}
  Consider the Kummer construction of a K3 surface~\( N \) as the
  resolution of \( X/\{\pm1\} \), where \( X=\Tn[4] = \bC^2/(\bZ+i\bZ)^2
  \), obtained by blowing up the \( 16 \) singular points.  Note that the
  Weierstrass \( \wp \)-function \( \wp\colon \bC\to\CP(1) \) descends to
  \( \Tn[2] \) and satisfies \( \wp(z)=\wp(-z) \), so each factor~\( \bC \)
  of \( \bC^2 \) defines a non-constant holomorphic map \( \wp_i \colon
  X/\{\pm1\}\to\CP(1) \), which we may then pull-back to the
  desingularisation~\( N \).  The classes \( \bigl[\Phi(\wp_i)\bigr] \) are
  inequivalent in \( \Horb^2(X/\{\pm1\},\bZ) = (H^2(X,\bZ))^{\pm1} \),
  indeed \( \bigl[\Phi(\wp_1)\wedge \Phi(\wp_2)\bigr] \) is non-zero in \(
  H^4(N) \).

  Taking \( F_i = \Phi(\wp_i) \), \( i=1,2 \), we see that \( M = N\times
  \Tn[2] \) may be twisted to an \SKT \( 6 \)-manifold \( W \) with finite
  fundamental group.  Taking the universal cover, we thus obtain a compact
  simply-connected \SKT manifold, that is a \( \Tn[2] \)-bundle over the K3
  surface \( N \).  We will see below that \( W \) is not Kähler.
\end{example}

These examples generalise as follows, cf. \cite{Grantcharov-GP:CY-toric}
and \cite{Goldstein-P:SU3}.

\begin{proposition}
  \label{prop:Kaehler-torus}
  Let \( N \) be a compact simply-connected \SKT manifold.  Suppose that
  for some even integer \( n=2k \), there are \( n \) closed integral \(
  (1,1) \)-forms \( F_i\in\Omega_\bZ^{1,1}(N) \) with \( [F_i]\in
  H^2(N,\bR) \) linearly independent and such that \( \sum_{i,j=1}^n
  \gamma_{ij}F_i\wedge F_j = 0 \) for some positive definite matrix \(
  (\gamma_{ij}) \in M_n(\mathbb R) \).  Then there is a compact
  simply-connected \( \Tn \)-bundle \( \widetilde W \) over~\( N \) whose
  total space is \SKT.

  Moreover,
  \begin{enumerate}
  \item no complex structure on \( \widetilde W \) compatible with the
    fibration to \( N \) is of Kähler type;
  \item if \( b_2(N)=n \), then the topological manifold \( \widetilde W \)
    admits no Kähler metric.
  \end{enumerate}
\end{proposition}

\begin{proof}
  We take \( W \) to be the twist of \( N\times \Tn \), where the flat
  Kähler metric on \( \Tn=\bR^n/\bZ^n \) is given by \( (\gamma_{ij}) \)
  with respect to the standard generators with a compatible complex
  structure: from the classification of quadratic forms \( (\gamma_{ij}) =
  Q^TQ \) and we may then take \( I = Q^{-1}I_0Q \), where \( I_0 \) is the
  standard complex structure on \( \bR^n=\bC^k \).  Then the discussion
  above shows that \( W \) is \SKT.  Topologically \( W \) is a principal
  torus bundle over \( N \) with Chern classes~\( [F_i] \).  As the \(
  [F_i] \) are linearly independent, the exact homotopy sequence shows that
  \( W \)~has finite fundamental group and its universal cover \(
  \widetilde W \) is a compact simply-connected \SKT manifold.
  Furthermore, \( \widetilde W \) is itself a \( \Tn \)-bundle over \( N \)
  with Chern classes \( [F_i'] \), rational linear combinations of the~\(
  [F_i] \), linearly independent over~\( \bR \).

  For the second part, the projection \( \widetilde W\to N \) is
  holomorphic, with complex fibres that are homologous to zero.  But any
  complex submanifold of a Kähler manifold is non-zero in homology, so the
  complex structure on \( \widetilde W \) can not be of Kähler type.

  For the last part, we use the Leray spectral sequence to compute \(
  b_2(\widetilde W) \).  We have \( E_2^{p,q} = H^p(N)\ot H^q(\Tn)
  \).  As the \( [F_i'] \) are linearly independent, we see that the map \(
  d_2\colon E_2^{0,q}\to E_2^{2,q-1} \), \(
  d_2[\theta_{i_1,\dots,i_q}']=\sum_{j=1}^q
  (-1)^{j+1}[F_{i_j}']\ot[\theta_{i_1,\dots,\hat{i_j},\dots,i_q}'] \) is
  injective for each \( q\geqslant 1 \).  In particular, \( E_3^{0,2} =
  \{0\} \) and \( d_3 \) is zero on \( E_3^{p,q} \) for \( p,q\leqslant 2
  \), so the spectral sequence stabilisers at level~\( E_3 \) for these
  terms.  This gives
  \begin{equation*}
    b_2(\widetilde W) = \sum_{i=0}^2\dim(E_3^{i,2-i}) = 0 + 0 + b_2(N)-n =
    b_2(N) - n. 
  \end{equation*}
  In particular, when \( b_2(N) = n \), we find \( b_2(\widetilde W) = 0 \)
  and so \( \widetilde W \) admits no Kähler metric.
\end{proof}

\begin{example}
  Let \( N_0 \) be a simply-connected projective Kähler manifold of real
  dimension~\( 4 \).  Fix an embedding \( N_0 \subset \CP(r) \).  Then a
  generic linear subspace \( \Proj(V) \subset \CP(r) \) of complex
  dimension \( r-2 \) meets \( N_0 \) transversely at a finite number of
  points \( p_1,\dots,p_d \).  Choose homogeneous coordinates \(
  [z_0,\dots,z_r] \) on \( \CP(r) \) so that \( \Proj(V) = (z_0=0=z_1) \).
  Then the blow-up  \( \widehat{\CP(r)} \) of \( \CP(r) \) along \(
  \Proj(V) \) is
  \begin{equation*}
    \widehat{\CP(r)} = \bigl\{ ([z_0,\dots,z_r],[w_0,w_1]) : z_0w_1 =
    z_1w_0 \bigr\} \subset \CP(r)\times\CP(1).
  \end{equation*}
  Let \( \pi_1 \) and \( \pi_2 \) denote the projections to the first and
  second factors of \( \CP(r)\times\CP(1) \).  Then \( N_1=\pi_1^{-1}(N_0)
  \) is the blow-up of~\( N_0 \) at \( p_1,\dots,p_d \), and \(
  f=\pi_2|_{N_1} \) is a non-constant holomorphic map from \( N_1 \) onto
  \( \CP(1) \).  As above, \( f \) defines \(
  F_1=\Phi(f)\in\Omega^{1,1}_\bZ(N_1) \), non-zero in cohomology with \(
  F_1^2=0 \).

  Iterating this construction, we find that for any \( n=2k \) there is a
  multiple blow-up \( N \) of~\( N_0 \) that satisfies the hypotheses of
  Proposition~\ref{prop:Kaehler-torus} in the form \( \sum_{i=1}^n F_i^2=0
  \), and hence \( N \) is the base of a simply-connected \SKT manifold of
  dimension \( n+4 \).
\end{example}

\begin{example}
  Let \( G \) be an even-dimensional compact Lie group.  This carries an
  \SKT structure with the complex structure of
  Samelson~\cite{Samelson:complex} and a bi-invariant metric
  \cite{Spindel-STvP:complex} .  Let \( T \) be the maximal torus of~\( G
  \).  Then one may twist \( G \) to \( W = G/T\times T \).  The space \(
  G/T \) is the maximal flag manifold and so Kähler, and \( T \)~has even
  rank.  By duality of the twist construction, we may thus obtain the \SKT
  structure of~\( G \) by twisting the Kähler manifold \( G/T\times T \).
\end{example}

\subsection{Non-instanton twists}
\label{sec:non-instanton}

If \( M \) is a torus, then \( M \) can be repeatedly twisted to produce
nilmanifolds.  Indeed every nilmanifold may be produced this way,
cf.~\cite{Swann:Bilbao}.  From the results of Fino, Parton \& Salamon
\cite{Fino-PS:SKT} for \SKT structures on \( 6 \)-dimensional nilmanifolds,
we can see that non-instanton twists are necessary to produce all such
examples.  We now demonstrate that there are simply-connected \SKT
manifolds obtainable from non\bdash instanton twists over K3 surfaces.

Let \( N \) be a K3 surface.  Yau's proof of the Calabi conjecture shows
that \( N \) admits hyperKähler metrics \( (g,\omega_I,\omega_J,\omega_K)
\).  From Looijenga's Torelli theorem \cite{Looijenga:Torelli} for period
maps of K3 surfaces there is a hyperKähler metric \( \tilde g \) whose
Kähler classes \( \tilde\omega_I \), \( \tilde\omega_J \), \(
\tilde\omega_K \) have integral periods.  Consider \( M=N\times \Tn[2] \)
as above, with the product Kähler metric for \( \Tn[2]=\bC/(\bZ+i\bZ) \)
standard and any Kähler metric on~\( N \) which is Hermitian with respect
to~\( I \) from the hyperKähler triple.  We twist \( M \) using \(
(F_1,F_2) = (F_1^0 + \tilde\omega_J, F_2^0 + \tilde\omega_K) \) with \(
F_i^0 \) of type \( (1,1) \) with respect to~\( I \).  Note that \(
\tilde\omega_J+i\tilde\omega_K \) is of type \( (2,0) \) with respect to~\(
I \), so the almost complex structure \( I_W \) on the resulting twist~\( W
\) is integrable by Lemma~\ref{lem:integrable}.  The condition for \( W \)
to be \SKT is now
\begin{equation*}
  d(IF_1)\wedge X_1^\flat + d(IF_2)\wedge X_2^\flat + F_1\wedge IF_1 +
  F_2\wedge IF_2 = 0.
\end{equation*}
We have \( (IF_1,IF_2) = (F_1^0 - \tilde\omega_J, F_2^0 - \tilde\omega_K)
\), so find \( d(IF_i)=0 \) and the \SKT condition reduces to \(
(F_1^0)^2+(F_2^0)^2 = \tilde\omega_J^2 + \tilde\omega_K^2 \).  As \(
\tilde\omega_I^2 = \tilde\omega_J^2 = \tilde\omega_K^2 \) for any
hyperKähler metric in dimension~\( 4 \), we can solve the \SKT equations by
taking \( F_i^0 = \varepsilon_i\tilde\omega_I \), \(
\varepsilon_i\in\{\pm1\} \).  The resulting \( \Tn[2] \)-twists are then
\SKT with finite fundamental group and their universal covers \( \widetilde
W \) are the promised manifolds.  As above the complex structure admits no
compatible Kähler metric.  One may compute \( b_2(\widetilde W) = 20 \) and
\( b_3(\widetilde W) = 42 \).  From this one can see that \( \widetilde W
\) is not a product of smaller dimensional manifolds.

\subsection{Complex volume forms}
\label{sec:volume}

If \( (M,I) \) is a complex manifold with trivial canonical bundle, it is
natural to ask which twists \( W \) also enjoy this property.

\begin{proposition}
  \label{prop:volume}
  Suppose the twist \( (W,I_W) \) of \( (M,I) \) via \( \xi \), \( F \)
  and~\( a \) is complex.  If \( (M,I) \) carries an invariant complex
  volume form~\( \Theta \), then the \Hrelated form \( \Theta_W \) is a
  complex volume form on \( (W,I_W) \) if and only if
  \begin{equation*}
    a^{-1}\xi\hook F^{1,1} = 0,
  \end{equation*}
  where \( F^{1,1}=\tfrac12(F+IF) \) is the \( (1,1) \)-part of~\( F \).
\end{proposition}

\begin{proof}
  Let \( \dim_{\bC} M = m \).  A section of \( \Lambda^{m,0} \) is
  holomorphic only if it is closed.  We have
  \begin{equation*}
    d\Theta_W \Hrel d\Theta - a^{-1} F\wedge\xi\hook\Theta =
    -a^{-1}F\wedge\xi\hook\Theta,
  \end{equation*}
  so we need to determine when the right-hand side vanishes.  This is a
  pointwise computation and we may thus use Lemma~\ref{lem:integrable}. 

  In the notation of Lemma~\ref{lem:integrable}, we wish to compute \(
  \sum_{i=1}^n F_i\wedge X_i\hook\Theta \).  Note that \( X\hook\Theta \)
  is of type \( (m-1,0) \) and that \( IX\hook\Theta = iX\hook\Theta \).
  Now, for \( k>2s \), we have \( F_k^{} = F_k^{(1,1)} \), whereas for \(
  j\leqslant s \), \( F_{2j-1}\wedge X_{2j-1}\hook \Theta + F_{2j}\wedge
  X_{2j}\hook \Theta = F_{(1)}\wedge X\hook\Theta + F_{(2)}\wedge IX\hook
  \Theta = F_{(c)}\wedge X\hook\Theta = F_{(c)}^{1,1}\wedge X\hook\Theta
  \), as \( F_{(c)}^{2,0}\wedge X\hook\Theta\in\Lambda^{m+1,0}=\{0\} \).
  However, for any \( (1,1) \)-form \( F \) we have \( F\wedge \Theta=0 \),
  so \( (X\hook F)\wedge \Theta + F\wedge (X\hook\Theta) = 0 \).  This gives
  \begin{equation*}
      \sum_{i=1}^n F_i\wedge X_i\hook\Theta
      = \sum_{i=1}^n F_i^{1,1}\wedge X_i\hook\Theta 
      = - \Bigl(\sum_{i=1}^nX_i\hook F_i^{1,1}\Bigr)\wedge \Theta.
  \end{equation*}
  The result follows from the fact that \( \cdot \wedge \Theta\colon
  \Lambda^1\to\Lambda^{m,1} \) is an \( \bR \)-linear isomorphism.
\end{proof}

\section{Hypercomplex and HKT geometry}
\label{sec:HKT}

We now turn geometries with multiple complex structures.  An \emph{almost
hypercomplex structure} on a manifold \( M \) is a triple of tangent bundle
endomorphisms \( I, J, K \in \End(TM) \) satisfying the identities
\begin{equation*}
  I^2=-1=J^2=K^2\qquad IJ = K = -JI.
\end{equation*}
The first identities say that \( I \), \( J \) and \( K \) are almost
complex structures.  If these three are integrable complex structures we
then have a \emph{hypercomplex structure} on~\( M \).  A metric \( g \)
satisfying \( g(IX,IY)=g(X,Y)=g(JX,JY)=g(KX,KY) \) for all \( X,Y\in TM \)
is said to be \emph{almost hyperHermitian}.

\begin{definition}
  An almost hyperHermitian structure \( (g,I,J,K) \) is \emph{hyperKähler
  with torsion} or \emph{\HKT} if
  \begin{equation}
    \label{eq:HKT}
    Id\omega_I = Jd\omega_J = Kd\omega_K,
  \end{equation}
  where \( \omega_I(X,Y) = g(IX,Y) \), etc.
\end{definition}

In \cite{Cabrera-S:aqH-torsion} it was shown that \HKT structures are
always hypercomplex, a condition that was previously included in the
definition of \HKT, cf.~\cite{Grantcharov-P:HKT}.  Equation~\eqref{eq:HKT}
is now the condition that the Bismut connections for \( (g,I) \), \( (g,J)
\) and \( (g,K) \) agree.

The results of \S\S\ref{sec:twist} and~\ref{sec:torsion} may now be applied
to these structures.  
Firstly, a direct consequence of
Proposition~\ref{prop:KT} and equation~\eqref{eq:HKT} is:

\begin{proposition}
  \label{prop:HKT}
  Suppose \( M \) is an \HKT manifold with twist data \( (\xi,F,a) \) and
  that the action \( \xi \) preserves the \HKT structure.  Then the twist
  \( W \) of \( M \) by \( (\xi,F,a) \) is \HKT if and only if
  \begin{equation}
    \label{eq:HKT-twist}
    a^{-1}IF\wedge \xi^\flat =
    a^{-1}JF\wedge \xi^\flat =
    a^{-1}KF\wedge \xi^\flat.
  \end{equation}
  \qed
\end{proposition}

\subsection{Instanton HKT twists}
\label{sec:instanton-HKT}

The condition~\eqref{eq:HKT-twist} is satisfied by any \emph{instanton},
meaning \( F\in S^2E\ot\aP \), where \( S^2E = \Lambda^{1,1}_I \cap
\Lambda^{1,1}_J \cap \Lambda^{1,1}_K \).  This immediately gives many
examples.

Consider the following building blocks; we need both integral Hamiltonian
instantons and torus symmetries.
\begin{asparaenum}
\item\label{item:tori} \emph{Tori} \( \Tn[4k] = \bH^k/\Lambda \), \(
  \bH=\bR^4 \), \( \Lambda\cong\bZ^{4k} \) a lattice.  These are
  hyperKähler, so \HKT.  They carry no Hamiltonian instantons, but supply
  symmetries.
\item\label{item:hK} \emph{Compact irreducible hyperKähler manifolds}.
  Passing to the universal cover, we may take these to be simply-connected.
  These give a rich supply of integral instantons.  One family of examples
  are provided by K3 surfaces~\( M \); here the instanton condition is just
  that \( F \) be self-dual; by the Torelli theorem there are examples
  where the integral instantons form a lattice of rank~\( 19 \).  Notice
  that these have no Killing vectors, since they are compact Ricci-flat and
  irreducible.
\item\label{item:group} \emph{Compact groups} and related homogeneous
  spaces.  If \( G \) is a compact, simple Lie group, then there is a torus
  \( \Tn[r] \) such that \( G\times \Tn[r] \) equipped with a bi-invariant
  metric is \HKT \cite{Grantcharov-P:HKT}.  The hypercomplex structure was
  determined by Joyce \cite{Joyce:hypercomplex} and the minimal value of~\(
  r \) may be found in \cite[Table~1]{Spindel-STvP:complex}.  In
  particular, \( \SU(2n+1) \) carries an \HKT structure for all~\( n \).
  These structures admit \HKT deformations with torus symmetry
  cf.~\cite{Pedersen-Poon:inhomogeneous}.  Also, some of these homogeneous
  \HKT structures descend to homogeneous spaces \( (G/H)\times \Tn[s] \),
  see~\cite{Opfermann-P:Hom-H+QKT}.  Again, the spaces \( G\times \Tn[r] \)
  carry no Hamiltonian instantons, but do supply symmetries.
\item\label{item:3S} \emph{Squashed \( 3 \)-Sasaki structures}.  A
  Riemannian manifold \( (\tS,h_0) \) of dimension \( 4n-1 \) is \( 3
  \)-Sasaki if the cone \( (\tS\times\bR, g_0 = dr^2+r^2h_0) \) is
  hyperKähler with complex structures invariant under~\( X =
  r\partial/\partial r \).  The metric \( h_0 \) is then Einstein with
  positive scalar curvature, so if \( \tS \) is compact then \( \pi_1(\tS)
  \) is finite.  Passing to the universal cover we may assume \( \tS \) is
  simply-connected.  Rescaling \( g_0 \) with different weights along and
  transverse to the quaternionic span of~\( X \), one may produce an \HKT
  metric \( g = dt^2+h \) that descends to \( \tS\times S^1 \) such that \(
  X \) acts a triholomorphic isometry \cite{Ornea-PS:potential}.  Galicki
  \& Salamon showed that harmonic \( 2 \)-forms for \( h_0 \) are
  instantons and orthogonal to the quaternionic span of~\( X \), so these
  are also instantons for the hypercomplex structure on \( \tS\times S^1
  \).  Thus any \( \tS \) with \( b_2(\tS)>0 \) provides integral
  Hamiltonian instantons on \( \tS\times S^1 \).  Many such examples of \(
  3 \)-Sasaki manifolds have been constructed by Boyer, Galicki and their
  coworkers, see~\cite{Boyer-G:book}; in particular, there are
  inhomogeneous examples in dimension~\( 7 \) with arbitrarily large \(
  b_2(\tS) \).  Moreover, such \( \tS \) often admit non-trivial isometries
  preserving both the \HKT structure of \( \tS\times S^1 \) and \( h_0 \)
  and hence the above mentioned instantons.  Thus the manifolds \(
  \tS\times S^1 \) provide a rich source of symmetries too, including free
  actions that have \( \xi\hook F \) non-zero.
\end{asparaenum}

The hypotheses of the following Theorem are satisfied by an \HKT space~\( M
\) that is a product of manifolds of the four types above and the
additional torus factors needed in types \itref{item:group}
and~\itref{item:3S}.  The existence of an appropriate Hamiltonian instanton
\( F \) is guaranteed if there are sufficiently many factors of type
\itref{item:hK} or type \itref{item:3S} with large \( b_2(\tS) \).

We say that a group action on a product \( A\times B \) \emph{projects
transitively} to~\( B \) if the action preserves the product structure, so
\( g\cdot(a,b)=(g\cdot a,g\cdot b) \), and the induced action on the second
factor~\( B \) is transitive.

\begin{theorem}
  \label{thm:HKT-examples}
  Let \( M = M_0 \times \Tn[p] \) be a compact \HKT manifold of dimension~\(
  4m \), \( p=b_1(M) \).  Suppose \( \Tn \) acts freely on \( M \)
  preserving the \HKT structure, projecting transitively to \( \Tn[p] \) and
  preserving a Hamiltonian instanton \( F\in\Omega^2_\bZ(M) \) of rank~\( n
  \) in~\( H^2(M,\tn) \).  Then there is a finite cover \( \widetilde W \)
  of a twist \( W \) of~\( M \) via \( F \) that is a compact
  simply-connected \HKT manifold.
\end{theorem}

Here the rank of \( F \) in~\( H^2(M,\tn) \) is the rank of \( [F] \) as a
linear map \( H^2(M,\bR) \to \tn^*\cong\bR^n \).

\begin{proof}
  From the discussion of the building blocks we see that the Hamiltonian
  instanton \( F \) is the pull-back of a form \( F_0 \) on~\( M_0 \).
  Then the twisting bundle \( P\to M = M_0\times \Tn[p] \) is \( P_0\times \Tn[p]
  \), where \( P_0\to M_0 \) is a principal \( \Tn \)-bundle with
  curvature~\( F_0 \).  As \( [F] \) has rank~\( n \), the fundamental
  group of \( P_0 \)~is finite.  The lifted \( \Tn \)-action \( \liftTM \)
  is free and projects transitively to \( \Tn[p] \).  In particular, \(
  W=P/\liftTM \) is diffeomorphic to a quotient of~\( P_0 \) by the free
  action of a compact group and so \( \pi_1(W) \) is finite.  However, by
  Proposition~\ref{prop:HKT}, the instanton condition also ensures that \(
  W \) is \HKT.
\end{proof}

For example, this Theorem gives \HKT metrics not only on torus bundles of
rank \( 4 \), \( 8 \), \( 12 \) and \( 16 \) over a single K3 surface, but
also on certain bundles with fibre \( \SU(n+1) \) or \( \tS\times S^1 \)
for a \( 3 \)-Sasaki manifold~\( \tS \).

Note that Verbitsky \cite{Verbitsky:hyperholomorphic-kaehler} constructs
\HKT metrics on vector bundles out of instanton connections, but the
resulting compact quotients are never simply-connected.

\subsection{Special Obata holonomy}
\label{sec:SLnH}

If \( (M^{4m},I,J,K) \) is any hypercomplex manifold, then there is a
unique torsion-free connection \( \Ob \) preserving the complex structures
\cite{Obata:connection}.  This implies that the holonomy of \( \Ob \) lies
in the group \( \GL(m,\bH) \).  Verbitsky~\cite{Verbitsky:holonomy} showed
using a version of Hodge theory for \HKT manifolds, that if \( M \) is \HKT
and the canonical bundle is trivial then the holonomy reduces to the
subgroup \( \SL(m,\bH) \).  The following is a simpler result that will be
sufficient for our purposes.

\begin{proposition}
  \label{prop:SLm}
  Let \( (M^{4m},I,J,K) \) be a hypercomplex manifold with trivial
  canonical bundle.  If \( M \) admits a complex volume form \( \Theta \)
  with respect to \( I \) that satisfies \( J\Theta = \overline \Theta \),
  then the Obata connection has holonomy in \( \SL(m,\bH) \).
\end{proposition}

\begin{proof}
  As \( \Ob \) preserves the complex structure \( I \), we have \(
  \Ob\Theta= \theta\otimes \Theta \) for some one-form \(
  \theta\in\Omega^1(M,\bC) \).  Now \( \Ob \)~is torsion-free, so \(
  d\Theta \) is the alternation of \( \Ob\Theta \) and \(
  0=d\Theta=\theta\wedge\Theta = \theta^{0,1}\wedge\Theta \in
  \Omega^{2m,1}(M) \), where type decompositions are with respect to~\( I
  \).  We conclude that \( \theta^{0,1}=0 \).

  Now \( \Ob \) preserves \( J \) too, so \( \theta(X)J\Theta =
  J\Ob_X\Theta = \Ob_X(J\Theta) = \Ob_X\overline\Theta =
  \overline{\theta}(X)J\Theta \) for any tangent vector~\( X \).  This
  shows that \( \theta \) is a real one-form, and from \( \theta^{0,1}=0 \)
  we find that \( \theta=0 \).  Thus \( \Theta \) is Obata parallel and \(
  \Theta\wedge\overline\Theta \) is a parallel real volume form for \( \Ob
  \).
\end{proof}

\begin{theorem}
  \label{thm:SLm-examples}
  The spaces of Theorem~\ref{thm:HKT-examples} built using factors of type
  \itref{item:tori} and~\itref{item:hK} are simply-connected hypercomplex
  \( 4m \)-manifolds with holonomy contained in \( \SL(m,\mathbb H) \).
\end{theorem}

\begin{proof}
  The given factors are hyperKähler, so \( \Theta = (\omega_J+i\omega_K)^m
  \) is a complex volume form on~\( M=M_0\times \Tn[p] \).  In our case \( M_0
  \) is a product of compact irreducible hyperKähler manifolds, the
  twisting form \( F \) is the pull-back of a form on~\( M_0 \), whereas
  the group action of \( \TM \) is simply that of the factor \( \Tn[p] \).  So
  \( a^{-1}\xi\hook F=0 \) and Proposition~\ref{prop:volume} implies that
  \( \widetilde W \) has holomorphically trivial canonical bundle with
  complex volume form \( \Theta_W = (\omega_J^W+i\omega_K^W)^m \).  But \(
  J\Theta_W = (J\omega_J^W+iJ\omega_K^W)^m = (\omega_J^W-i\omega_K^W)^m =
  \overline{\Theta_W} \).  The result now follows from
  Proposition~\ref{prop:SLm}.
\end{proof}

In particular, this constructs examples on instanton torus bundles over
products of K3 surfaces.

The other building blocks of type \itref{item:group} and~\itref{item:3S}
described above typically do not have \( c_1=0 \) and so can not be used in
the construction of Theorem~\ref{thm:SLm-examples}.

\subsection{Non-instanton HKT twists}
\label{sec:non-instanton-HKT}

Let us note that \HKT twists of the type in Proposition~\ref{prop:HKT} need
not come from instantons, at least in the non-compact case.

\begin{example}
  Let \( M = \bR_{>0}\times \Tn[3] \subset \bH/(\bZ i+\bZ j+\bZ k) \) with the
  flat hyperKähler structure.  Let \( X_0 = \partial/\partial x^0 \) be the
  generator of the first factor, so that \( X_1=IX_0 \), \( X_2=JX_0 \) and
  \( X_3=KX_0 \) generated the three circle factors of \( \Tn[3] \).  Let \(
  b_0 \), \( b_1 \), \( b_2 \), \( b_4 \) be the (unit length) dual
  one-forms.  Put \( F_I = b_0\wedge b_1 = b_{01} \), \( F_J = b_{02} \)
  and \( F_K = b_{03} \).  Then \( X_1\hook F_I = -b_0 = -dx^0 = X_2\hook
  F_J = X_3\hook F_K \), with all other \( X_i\hook F_A \) zero.  We may
  thus take \( a = -\diag(x^0,x^0,x^0) \).  Then
  \begin{equation*}
    \begin{split}
      a^{-1}IF\wedge \xi^\flat
      &= -\tfrac1{x^0}(IF_I\wedge b_1+IF_J\wedge b_2 +IF_K\wedge b_3)\\
      &= -\tfrac1{x^0}(0+b_{13}\wedge b_2 -b_{12}\wedge b_3)\\
      &= \tfrac2{x^0}b_{123} = a^{-1}JF\wedge\xi^\flat =
      a^{-1}KF\wedge\xi^\flat. 
    \end{split}
\end{equation*}
Thus we may twist to obtain a new \HKT metric, even though \( F\notin
S^2E\ot\aP \) since \( JF_I\ne F_I \).
\end{example}

\subsection{Hypercomplex manifolds that are not HKT}
\label{sec:hc-not-HKT}

Proposition~\ref{prop:N} directly tells when hypercomplex structures are
preserved by the twist construction.

\begin{proposition}
  \label{prop:hc}
  Suppose that \( M \) is hypercomplex with twist data \( (\xi,F,a) \) and
  that the action \( \xi \) preserves the hypercomplex structure.  Then the
  twist \( W \) is hypercomplex if and only if
  \begin{equation*}
    \cL_I\cF = \cF = \cL_J\cF = \cL_K\cF.
  \end{equation*}
  \qed
\end{proposition}

Let us use this to give examples of compact simply-connected hypercomplex
manifolds that are not \HKT.

Let \( N \) be a hyperKähler K3 surface \( N \) for which the three Kähler
forms \( \omega_I \), \( \omega_J \) and \( \omega_K \) have integral
periods, cf.~\ref{sec:non-instanton}.  Take \( M \) to be the product \( N
\times \Tn[4] \) with the hyperKähler torus \( \Tn[4] = (S^1)^4 = \bH/\bZ^4
\) and let \( X_0 \), \( X_1=IX_0 \), \( X_2=JX_0 \), \( X_3=KX_0 \) be
vector fields generating the four circles.  Let \( \omega_0 \) be any
non-zero self-dual element of \( \Omega^2_\bZ(N) \).  Then \( \omega_0 \)
is of type \( (1,1) \) with respect to each complex structure.

Now twist~\( M \) by \( \cF = X_0\ot F_0 + X_1\ot F_I + X_2\ot F_J + X_3\ot
F_K \), where \( F_A = \pi_N^*\omega_A \) is the pull-back of~\( \omega_A
\).  Since \( X_i\hook F_A=0 \) we may take the twisting function~\( a \)
to be the identity matrix.  The resulting twist \( W \) has a finite
fundamental group and so its universal cover~\( \widetilde W \) is
simply-connected.  We have
\begin{equation*}
  \begin{split}
    \cF &= X_0 \ot F_0 + X_1 \ot F_I + \re((X_2-iX_3) \ot
    (F_J+iF_K))\\ 
    &\in T\ot\Lambda^{1,1}_I + \rcomp{T_I^{1,0}\ot \Lambda_I^{2,0}},
  \end{split}
\end{equation*}
so \( \cL_I\cF=\cF \), and similarly for \( I \), \( J
\) and \( K \).  Thus, \( W \) is hypercomplex, by Proposition~\ref{prop:hc}.
However, Proposition~\ref{prop:HKT} shows that the geometry on~\( W \) is not
\HKT, since \( M \) is \HKT and
\begin{equation}
  \label{eq:notHKT}
  \begin{split}
    &a^{-1}IF\wedge \xi^\flat = F_0\wedge b_0 + F_I\wedge b_1 - F_J\wedge
    b_2 - F_K\wedge b_3\\ &\ne a^{-1}JF\wedge \xi^\flat = F_0\wedge b_0 -
    F_I\wedge b_1 + F_J\wedge b_2 - F_K\wedge b_3,
  \end{split}
\end{equation}
where \( b_i= X_i^\flat \).

I claim that the hypercomplex structure on \( \widetilde W \) admits no
compatible \HKT metric.  Suppose for a contradiction that \( g_0 \) is an
HKT metric on \( (\widetilde W,I,J,K) \).  As the hypercomplex structure on
\( W \) is constructed via the twist construction it is \( \Tn[4]
\)-invariant.  Averaging~\( g_0 \) over \( \pi_1(W) \) and then the
principal \( \Tn[4] \)-action gives a metric~\( g_1 \) on~\( W \) that is
still \HKT since the \HKT condition is linear and invariant under
tri-holomorphic pull-backs.  Untwisting \( W \) gives the original product
hypercomplex structure on \( N\times \Tn[4] \) and a hyperHermitian metric
\( g_2 \) on \( N\times\Tn[4] \) which is \( \Tn[4] \)-invariant.  As all
hyperHermitian metrics on a four-dimensional vector space are proportional,
we may write
\begin{equation*}
  g_2 = f\,g^{}_N + h \sum_{i=0}^3b_i^2 + 2\sum_{i=0}^3b_i\vee
  \alpha_i, 
\end{equation*}
where \( g_N \) is the hyperKähler metric on~\( N \), \(
\alpha_i\in\Omega^1(N) \), \( \alpha_1=I\alpha_0 \), etc., and \( f,h\in
C^\infty(N) \) are positive (pull-backs signs have be omitted).

Let \( \omega_A^{(i)} \) denote the Kähler forms of~\( g_i \).  Then
\begin{equation}
  \label{eq:HKTg1}
  \begin{split}
    I_Wd\omega_I^{(1)} &\Hrel Idf\wedge \omega_I + Idh \wedge
    (b_{01}+b_{23})\eqbreak - 2\sum_{i=0}^3 b_i\wedge Id\alpha_i +
    \sum_{i=0}^3IF_i\wedge(h\, b_i + \alpha_i).
  \end{split}
\end{equation}
For \( W \) \HKT, we have \( I_Wd\omega_I^{(1)} = J_Wd\omega_J^{(1)} =
K_Wd\omega_K^{(1)} \).  Using~\eqref{eq:HKTg1} and considering the
coefficient of \( b_{01} \) gives \( Idh = 0 \), so \( h \) is constant.
Looking at the coefficient of \( b_i \), we find
\begin{equation}
  \label{eq:dalpha}
  \begin{split}
    2(Id\alpha_i) - h\, IF_i &= 2(Jd\alpha_i) - h\, JF_i \\
    &= 2(Kd\alpha_i) - h\, KF_i.
  \end{split}
\end{equation}
For \( i=0 \), this gives \( Id\alpha_0 = Jd\alpha_0=Kd\alpha_0 \), so \(
d\alpha_0 \) is self-dual on the compact space~\( N \) and therefore zero.
As \( b_1(N) = 0 \), we may write \( \alpha_0 = d\phi \), for some \(
\phi\in C^\infty(N) \).  

Equation~\eqref{eq:dalpha} for \( i=1 \), is
\begin{equation}
  \label{eq:alpha1}
  2(Id\alpha_1) - h\omega_I = 2(Jd\alpha_1) + h\omega_I = 2(Kd\alpha_1) + h\omega_I.
\end{equation}
This first gives that \( Jd\alpha_1 = Kd\alpha_1 \), so \( d\alpha_1 \in
\Lambda^{1,1}_I \).  Writing \( d\alpha_1 = \beta + \lambda\omega_I \),
with \( \beta\in\Omega_+^2(N) \), gives \( Jd\alpha_1= \beta -
\lambda\omega_I \) and \eqref{eq:alpha1} implies \( \lambda = h/2 \),
which is constant.

However, \( \alpha_1 = I\alpha_0 = Id\phi \) and so \( \lambda =
\Delta_N\phi \), where \( \Delta_N \) is the Laplacian of the hyperKähler
metric~\( g_N \).  Since \( \Delta_N \) has image orthogonal to the
constant functions, we conclude that \( \lambda=0 \) and hence \( h=0 \),
contradicting the positive definiteness of~\( g_2 \).  Thus \( g_1 \) can
not be \HKT and the hypercomplex structure on~\( \widetilde W \) admits no
compatible \HKT metric.

To summarise:

\begin{theorem}
  There are compact simply-connected hypercomplex \( 8 \)-manifolds that
  admit no compatible \( \HKT \) metric.  Moreover these exist with Obata
  holonomy contained in \( \SL(2,\bH) \).
\end{theorem}

\begin{proof}
  It remains to prove the final assertion.  However, the twist \( W \)
  constructed above starts from an hyperKähler manifold and has \(
  a^{-1}\xi\hook F^{1,1} = 0 \), so as in the proof of
  Theorem~\ref{thm:SLm-examples} the Obata holonomy reduces.
\end{proof}

\begin{theorem}
  There are non-trivial compact simply-connected hypercomplex manifolds in
  all dimensions \( 4m\geqslant 8 \) that admit no compatible \( \HKT \)
  metric.  Furthermore, examples exist with holonomy in \( \SL(m,\bH) \).
\end{theorem}

The examples constructed are torus bundles over a product base and by
`non-trivial' we mean that the structure does not split as a product of
torus bundles over a product of factors of the base.

\begin{proof}
  Let \( W_0 \) be the \( 8 \)-dimensional example constructed above.
  Untwist the symmetry \( X_0 \) by \( F_0 \) to get the hypercomplex
  manifold \( W_1 = S^1\times B_7 \).  Let \( M_2 \) be a product of \(
  (m-2) \) factors that are each K3 surfaces with a fixed choice of
  non-zero integral self-dual \( 2 \)-form \( F_{(i)} \).  Twisting the \(
  S^1 \) factor of \( W_1 \) by \( F = F_0+F_{(1)}+\dots+F_{(m-2)} \) we a
  obtain a \( 4m \)-dimensional hypercomplex manifold \( W_2 \) with finite
  fundamental group and Obata holonomy in~\( \SL(m,\bH) \).  The universal
  cover~ \( \widetilde{W_2} \) of \( W_2 \) is the required example.

  The hypercomplex manifold \( \widetilde{W_2} \) admits no compatible \HKT
  metric, since any potential \HKT metric may be averaged so that it
  descends to a torus invariant \HKT metric on \( W_2 \) and then twisted
  to an \HKT metric on \( M\times W_0 \).  As \( W_0 \) is a hypercomplex
  submanifold of \( M\times W_0 \) such an \HKT metric would restrict to an
  \HKT metric on \( W_0 \) itself, but that is a contradiction.  Thus \(
  \widetilde{W_2} \) carries no compatible \HKT metric.
\end{proof}

Using the squashed \( 3 \)-Sasaki building blocks of
\S\ref{sec:instanton-HKT} one may obtain examples with \( c_1 \) non-zero
and so Obata holonomy not in \( \SL(m,\bH) \).

In dimension~\( 8 \) we can prove a more general non-triviality result.

\begin{theorem}
  Any compact hypercomplex \( 8 \)-manifold that does not admit an \HKT
  metric is not a non-trivial product of smaller dimensional hypercomplex
  manifolds.
\end{theorem}

\begin{proof}
  Suppose that \( M^8 \) is such a product \( N_1\times N_2 \) with \( N_i
  \) hypercomplex.  Then each \( N_i \) has dimension~\( 4 \).  However
  Boyer~\cite{Boyer:hH} showed that any compact hypercomplex \( 4
  \)-manifold is either \( \Tn[4] \), a K3 surface or \( S^3\times S^1 \).
  But each of these examples admits a compatible \HKT metric, and the
  product structure on \( M = N_1\times N_2 \) is then \HKT, contradicting
  the hypotheses.
\end{proof}

\providecommand{\bysame}{\leavevmode\hbox to3em{\hrulefill}\thinspace}
\providecommand{\MR}{\relax\ifhmode\unskip\space\fi MR }
\providecommand{\MRhref}[2]{%
  \href{http://www.ams.org/mathscinet-getitem?mr=#1}{#2}
}
\providecommand{\href}[2]{#2}

\end{document}